\documentclass[12pt,a4paper]{amsart}

\usepackage[latin1]{inputenc}
\usepackage{amsfonts,amssymb}
\usepackage{amstext}
\usepackage{amsthm}
\usepackage{amsmath}
\usepackage{amscd}
\newtheorem{thm}{Theorem}
\newtheorem{lem}[thm]{Lemma}
\theoremstyle{definition}

\newtheorem{ex}{exercice}
\newcommand{\be}{\begin{ex} \normal }\newcommand{\ee}{\end{ex}}

\numberwithin{equation}{section}

\newcommand{\N}{\mathbb N}
\newcommand{\R}{\mathbb R}

\newenvironment{appli}{\left( \begin{array}{ccc}}{\end{array} \right)}

\newcommand{\ba}{\begin{appli}}
\newcommand{\ea}{\end{appli}}

\newcommand{\ds}{\displaystyle}

\title{Compact metrizable groups are isometry groups of compact metric spaces}
\author{Julien Melleray}

\begin{document}
\begin{abstract}
This note is devoted to proving the following result: given a
compact metrizable group $G$, there is a compact metric space $K$
such that $G$ is isomorphic (as a topological group) to the
isometry group of $K$.
\end{abstract}
\maketitle
$ $\\

\noindent {\bf Introduction } \\

\noindent It is an easy fact that, if $X$ is a Polish metric space
and $Iso(X)$ is endowed with the topology induced by the product
topology on $X^X$, then  it turns $Iso(X)$ into a Polish group.
One may then wonder whether all Polish groups are of that form;
this was proved by Gao and Kechris in \cite{gaokec}:

\begin{thm} (Gao-Kechris)\label{gaokechris}
Let $G$ be a Polish group. Then there exists a Polish metric space
$Y$ such that $G$ is isomorphic to $Iso(Y)$.
\end{thm}

\noindent Similarly, it is easy to see that, if $K$ is a compact
metric space, then $Iso(K)$ is a compact metrizable group. \\
Given what we saw above, it is natural to wonder whether the
converse holds; this question was mentioned to me by Alekos
Kechris (private correspondence). The aim of this note is to
provide a positive answer to that question:

\begin{thm} \label{main}
Let $G$ be a compact metrizable group. Then there exists a compact
metric space $K$ such that $G$ is isomorphic to $Iso(K)$.
\end{thm}

\noindent To obtain the proof, it is natural to  try to find a
simpler proof of theorem \ref{gaokechris} than the one in
\cite{gaokec}, which is a little buried among other considerations
(it is a byproduct of a proof,
so no "direct" proof is given).\\
It turns out that such a simpler proof exists, and it is not very
hard to use a variation of it in order to prove theorem \ref{main}. \\
The paper is organized as follows: first we give a simple proof of
theorem \ref{gaokechris}, then we show how to adapt this proof in
order to obtain theorem \ref{main}.\\

\noindent \textit{Acknowledgements}: I would like to thank Alekos
Kechris, who told me about the problem studied here.\\
$ $ \\

\noindent {\bf Notations and Definitions}\\

\noindent If $(X,d)$ is a complete separable metric space, we will
say that it is a \textit{Polish
metric space}, and will often write it simply $X$. \\
A \textit{Polish group} is a topological group whose topology is
Polish; if $X$ is a separable metric space, then we will denote
its isometry group by $Iso(X)$, and endow it with the pointwise
convergence topology, which turns it into a second countable
topological group, and into a Polish group if $X$ is Polish
(see \cite{beckec} for a thorough introduction to the theory of Polish groups). \\

\noindent Let $(X,d)$ be a metric space; we will say that $f: X
\to \R $ is a \textit{Kat\v{e}tov map } iff
$$\forall x,y \in X \ |f(x)-f(y)|\leq d(x,y) \leq f(x)+f(y)\ \ .$$
These maps correspond to one-point metric extensions of $X$. We
denote by $E(X)$ the set of all Kat\v{e}tov maps on $X$; we endow
it with the sup-metric, which turns it into a complete
metric space.\\
If $f \in E(X)$ and $S \subseteq X$ are such that
$f(x)=\inf(f(s)+d(x,s) \colon s \in S)$, we say that $S$ is a
support of $f$; or that $S$ \textit{controls} $f$ .\\
It is useful to note here the following easy fact: if $f,g \in
E(X,\omega)$ are both supported by some set $S$, then
$d(f,g)=\sup_{s\in S} |f(s)-g(s)|$.
\\
Notice that $X$ isometrically embeds in $E(X)$ via the Kuratowski
embedding $x \mapsto \delta_x$, where $\delta_x(y)=d(x,y)$, and
that one has, for any $f\in E(X)$, that
$d(f,\delta_x)=f(x)$.\\
$ $ \\
\noindent {\bf Proofs of the theorems: } \\

\noindent As promised, we begin by proving theorem
\ref{gaokechris}. \\
Let $G$ be a Polish group, and $d$ be a left-invariant distance on
$G$. Let $X$ be the completion of $(G,d)$. Then the left
translation action of$G$ on itself extends to an action by
isometries of $G$ on $X$, and it is not hard to chaeck that this
provides a continuous embedding of $G$ in $Iso(X)$. We identify
$G$ with the corresponding (closed) subgroup of $Iso(X)$, and make
the additional assumption that $X$ is of diameter $\leq 1$.\\

\noindent {\bf Claim.} For all $\varphi\in Iso(X) \setminus G$,
there exist $x_1,\ldots x_m \in X $ and $\varepsilon >0$ such that
$V_{\varphi}=\{\psi \in Iso(X) \colon \forall 1\leq k \leq m \
d(\psi(x_{k}),\varphi(x_{k})) < \varepsilon \} \subseteq Iso(X)
\setminus G$, and $2m \varepsilon = \min (d(x_{l},x_{k}))$ . \\

\noindent {\bf Proof.} Obvious.\\

\noindent Choose for all $\varphi \in Iso(X) \setminus G$ such a
$V_{\varphi}$; then, since $Iso(X)$ is Lindelöff, there are
$\{\varphi_i\}_{i \geq 1}$ such that $Iso(X) \setminus G=
\bigcup_{i \geq 1} V_{\varphi_i}$. We denote $V_{\varphi_i}=\{\psi
\in Iso(X) \colon \forall 1\leq k \leq m_i \ d(y^i_k,\psi (x^i_k))
< \varepsilon_i \}$.\\

\noindent For each $i \geq 1$ we define maps $f_i,g_i \in E(X)$:\\
$f_i(x)= \min \big ( \min_{1 \leq k \leq m_i}  (1+ d(x,x^i_k)  +
2(k-1) \varepsilon_i), 1+ 2m_i \varepsilon_i \big ), \ \mbox{ and
}$\\
$g_i(x)= \min \big ( \min_{1 \leq k \leq m_i}  (1 + d(x,y^i_k)+
2(k-1) \varepsilon_i )\, , 1 +2m_i \varepsilon_i \big )\ .\qquad \
\ \ $

\noindent If $ \varphi \in Iso(X)$, we let $\varphi^*$ denote its
(unique) extension to $E(X)$, defined by
$\varphi^*(f)(x)=f(\varphi^{-1}(x))$. We now have the following
lemma :

\begin{lem} \label{lem1}
$ \ds{\forall \varphi \in Iso(X)\, \forall i \geq 1 \ \big
(\varphi \in V_{\varphi_i} \big ) \Leftrightarrow \big
(d(\varphi^{*}(f_i),g_i) < \varepsilon_i \big )\ .}$
\end{lem}

\noindent {\bf Proof of Lemma \ref{lem1}.} \\
Let $ \varphi \in V_{\varphi_i}$. Then the various inequalities
involved imply that:
$$\varphi^{*}(f_i)(y^i_k) \ = 1+  d(\varphi(x^i_k),
y^i_k)+ 2(k-1) \varepsilon_i ,\mbox{ so }
|\varphi^{*}(f_i)(y^i_k)- g_i(y^i_k)| \ \ \  < \varepsilon_i \,
.$$
$$g_i((\varphi(x^i_k))= 1 + d(\varphi(x^i_k),
y^i_k) + 2(k-1) \varepsilon_i, \mbox{ so }
|\varphi^{*}(f_i)(\varphi(x^i_k))- g_i(\varphi(x^i_k))| <
\varepsilon_i \, .$$

\noindent Since $\varphi^*(f_i)$ and $g_i$ are both supported by
the set $\{\varphi(x^i_k)\}_{k=1\ldots m_i} \cup
\{y^i_k\}_{k=1\ldots m_i}$, the inequalities above are enough to
show that $d(\varphi^*(f_i),g_i) < \varepsilon_i$. \\
Conversely, let $\varphi \in Iso(X)$ be such that
$d(\varphi^{*}(f_i),g_i) < \varepsilon_i$. We prove by induction
on $k=1 \ldots m$ that $d(\varphi(x^i_k),y^i_k) < \varepsilon_i$. \\
To see that this is true for $k=1$, remark that we have that
$g_i(\varphi(x^i_1)) < 1 + \varepsilon_i$, so that we must have
$g_i(\varphi (x^i_1)) = g_i(y^i_1) + d(y^i_1,x^i_1) $ .
In turn, this implies that  $ d(y^i_1,x^i_1)<  \varepsilon_i $. \\
Suppose now that we have proved the result up to rank $k-1 \leq
m-1$. \\
Notice that we have this time that  $g_i(\varphi(x^i_k)) < 1
 + \varepsilon_i +2(k-1) \varepsilon_i$ (*). Also,
we know that for all $l < k$ we have $d(x^i_k,x^i_l) >
2m_i\varepsilon_i$, and $d(x^i_k,x^i_l) < \varepsilon_i$. Thus,
$d(x^i_k, y^i_l) > (2m_i -1) \varepsilon_i$ for all $l< k$. \\
It is then clear that (*) implies that $g_i(\varphi(x^i_k)) =
g_i(y^i_k) + d(y^i_k,x^i_k) $, which means that
$d(y^i_k,x^i_k)< \varepsilon_i \, . \hfill \lozenge$ \\

\noindent We now let $F_0=X$, $F_i=\overline{\{\varphi^{*}.f_i
\colon g \in G\}} \subset E(X)$, and $\ds{Z=\overline{\cup F_i}}
$. Notice that $Z$ is a Polish metric space, since it is closed in
$E(X)$, which is complete, and it admits a separable dense subset.
\\
Also, it is important to remark that lemma \ref{lem1} implies that
$d(\psi^*(f_i),g_i) \geq \varepsilon_i$ for all $\psi \in G$; in
other words, $d(F_i,g_i) \geq
\varepsilon_i$ for all $i$. \\

\begin{lem}\label{lemme2}
Any element $ \varphi$ of $G$ extends (uniquely) to an isometry
$\varphi^Z$ of $Z$, and
$$\{\varphi^Z \colon \varphi \in G \} = \{ \varphi \in Iso(Z) \colon \varphi(X)=X \mbox{ and }\forall i \geq 0 \ \varphi(F_i)=F_i\}$$
\end{lem}

\noindent {\bf Proof of Lemma \ref{lemme2} .} \\
The first assertion is easy to prove: since $ \varphi^*(F_i)=F_i$
for all $\varphi \in G$, we see that $\varphi^*(Z)=Z $ for all
$\varphi \in G$. The fact that this extension is unique is a
classical consequence of the definition of the distance on $E(X)$,
see for
instance \cite{Katetov} .\\
It is also clear that $\{\varphi^Z \colon \varphi \in G \}
\subseteq \{ \varphi \in Iso(Z) \colon \varphi(X)=X \}$; to prove
the converse, take $\varphi \in Iso(Z)$ such that $\varphi(X)=X$ and $\varphi(F_i)=F_i$ for all $i$. \\
So $\varphi_{|_X}$ is an isometry of $X$ such that
$d\big((\varphi_{|_X})^*(f_i),g_i\big)=d(\varphi(f_i),g_i) \geq
\varepsilon_i$ for all $i$, which means that $\varphi_{|_X} \not
\in V_{\varphi_i}$ for all $i$, so that $\varphi_{|_X} \in G$ and
we are done.
$\hfill \lozenge$\\

\noindent To conclude the proof of the theorem, notice that $Z$ is
a bounded Polish metric space, so we may assume that
$\mbox{diam(Z)} \leq 1$.\\
Then we conclude as in \cite{gaokec}: for each $i \in I$ we let
$y_i \in E(Z)$ be defined by $y_i(z)=d(y_i,z)=(i+2)+d(z,F_i)$
for all $z \in Z$. Notice that $\varphi^*(y_i)=y_i$ for all $i \in I$ and all $\varphi \in G$.\\
Let $Y= Z \cup \{y_i\} \subset E(Z)$; $Y$ is complete, and we
claim that $G$ is isomorphic to $Iso(Y)$. \\
Indeed, any element $\varphi $ of $G$ has a unique isometric
extension $\varphi^Y$ to $Y$, and the mapping $\varphi \mapsto
\varphi^Y$ is continuous. Conversely, let $\psi$ be an isometry of
$Y$; necessarily $\psi(Z)=Z$, and $\psi(y_i)=y_i$ for all $i$. \\
Since $F_i=\{z \in Z \colon d(z,y_i)= \}$, we must have
$\psi(F_i)=F_i$ for all $i$. So, we see that there is some
$\varphi \in G$ such that $\psi_{|_Z}= \varphi^Z$, so that $\psi = \varphi^Y $. \\
To conclude the proof of theorem \ref{gaokechris}, recall that any
bijective, continuous morphism between two Polish groups is
actually bicontinuous.$\hfill
\lozenge$ \\

\noindent Now we will see that the ideas of this proof enable one to also prove theorem \ref{main}:\\

\noindent {\bf Proof of theorem \ref{main}.} \\
We may again assume that $G$ has more than two elements. Let $d$
be an invariant metric on $G$: the metric space $X=(G,d)$ is
compact, and $G$ embeds topologically in $Iso(X)$, via the mapping
$g \mapsto (x \mapsto g.x)$. We again identify $G$ with the
corresponding (closed) subgroup of $Iso(X)$, and make
the additional assumption that $X$ is of diameter $\leq 1$. \\

\noindent We again choose for, all $\varphi \in Iso(X) \setminus
G$, a $V_{\varphi}$ as in the claim; there are $\{\varphi_i\}_{i
\geq 1}$ such that $Iso(X) \setminus G= \bigcup_{i \geq 1}
V_{\varphi_i}$. We again denote $V_{\varphi_i}=\{\psi \in Iso(X)
\colon \forall 1\leq k \leq m_i \
d(y^i_k,\psi (x^i_k)) < \varepsilon_i \}$.\\
For each $i \geq 1$ we define slightly different maps $f_i,\ g_i$:
$$f_i(x)= \min \big ( \min_{1 \leq k \leq m_i}  (1+\frac{1}{2^i}+ d(x,x^i_k)  + 2(k-1) \varepsilon_i),
1+\frac{1}{2^i}+ 2m_i \varepsilon_i \big ), \ \mbox{ and }$$
$$g_i(x)= \min \big ( \min_{1 \leq k \leq m_i}  (1+ \frac{1}{2^i}\ + d(x,y^i_k) \  + 2(k-1) \varepsilon_i
)\, ,\  1+\frac{1}{2^i} +2m_i \varepsilon_i \big )\ .\qquad \ \ \
$$

\noindent If $ \varphi \in Iso(X)$, we let $\varphi^*$ denote its
(unique) extension to $E(X)$; we have again that
$$ \forall \varphi \in Iso(X)\, \forall i \geq 1 \ \big (\varphi
\in V_{\varphi_i} \big ) \Leftrightarrow \big
(d(\varphi^{*}(f_i),g_i) < \varepsilon_i \big )\ .$$

\noindent Now, we let $Y$ be the set of $f \in E(X)$ such that
$$\exists x_1 \ldots
x_n \ \ \forall x \ f(x)=\min \big( \min_{1\leq i\leq
n}(1+2(i-1)\varepsilon +d(x,x_i)), 1+2n \varepsilon\big),
$$ where $2n\varepsilon= \min(d(x_i,x_j) )$.\\
(for $n =0$ one gets $g$ defined by $g(x)=1$ for all $x \in X$.)

\begin{lem} $Y$ is compact.
\end{lem}

\noindent {\bf Proof.} Let $(f_i)$ be a sequence of maps in $Y$,
and let $x^i_1,\ldots x^i_{n_i}$ be points witnessing the fact
that $f_i \in Y$. \\
Then, either we can extract a sequence $f_{\varphi(i)}$ such that
$n_{\varphi(i)} \to + \infty$, or $(n_i)$ is bounded. \\
In the first case, notice that, since $X$ is totally bounded, one
must necessarily have that $\ds{\min_{1 \leq j < k \leq
n_{\varphi(i)}} d(x^{\varphi(i)}_j,x^{\varphi(i)}_k) \to 0}$ when
$n \to + \infty$, so that the definition of $f_i$ ensures that
$f_{\varphi(i)} \to g$.
\\
In the other case, we may extract a subsequence $f_{\psi(i)}$ such
that $n_{\psi(i)}=n$ for all $i$. \\
We assume that  $\min_{1 \leq j <k \leq n}
d(x^{\psi(i)}_j,x^{\psi(i)}_k) \geq \delta $ for some $\delta >0$
(if not, we can conclude as in the first case that some
subsequence of $(f_{\psi(i)})$ converges to $g$). But then, up to
another extraction, we can suppose that \\
$x^{\psi(i)}_1 \to x_1, \ldots,
 x^{\psi(i)}_n \to x_n$.\\
This implies that $\ds{\min_{1 \leq j <k \leq n}
d(x^{\psi(i)}_j,x^{\psi(i)}_k) \to \min_{1 \leq j <k \leq n}
d(x_j,x_k)}$, and one checks easily that
$f_{\psi(i)} \to f$ for some $f \in Y.$\\

\noindent We let, for all $i \geq 1$, $F_i=G^*.\{f_i\}$ (which is
a compact subset of
$E(X)$), and $Z=X \cup Y \cup \bigcup F_i$. \\
We proceed to prove  that $Z$ is compact: for that, it is enough
to see that any sequence $(x_n)$ of elements of $\cup F_i$ admits
a subsequence converging to some $z \in Z$. \\
We know by definition that $x_n=\varphi_n (f_{i_n})$ for some
$\varphi_n \in G$ and $i_n \in \N$. Since $G$ is compact, we may
assume that $\varphi_n \to \varphi$, so that it is enough to show
that $f_{i_n}$ admits a subsequence converging to some $z' \in Z$.
\\
We may of course assume that $i_n \to + \infty$. Notice that, by
definition, $f_{i_n}=\frac{1}{2^{i_n}} + h_n$, for some $h_n \in
Y$; since $Y$ is compact, and $i_n \to + \infty$, we
are done.\\

\noindent The end of the proof is very similar to that of theorem
\ref{gaokechris}, only a bit simpler (that is why we have chosen
the $f_i$ more carefully this time) : we pick $k \in E(Z)$ such
that $d(k,z)=\mbox{diam}(Z)+d(z,X)$, and let $\{k\} \cup Z=K$. \\
$K$ is compact, any element of $G$ extends uniquely to an isometry
of $K$, and the extension morphism is continuous.\\
So, we only need to prove that all isometries of $K$ are
extensions of elements of $G$; to that end, let $\psi \in
Iso(K)$.\\
We see that $\psi(k)=k,$ so $\psi(Z)=Z$. \\
Also, since $X=\{z \in Z \colon d(z,k)=\mbox{diam}(Z)\}$, we must
have $\psi(X)=X$. Similarly, $F_i=\{z \in Z \colon
d(z,X)=1+\frac{1}{2^i}\}$,
so $\psi(F_i)=F_i$.\\
We may now conclude as above: $\psi_{|_X} \not \in V_{\varphi_i}$
for all $i$, so $\psi_{|_X} \in G$, and we are done .$\hfill
\lozenge$.

$ $\\

\noindent  Equipe d'Analyse Fonctionelle, Universit\'e Paris 6 \\
Bo\^ite 186, 4 Place Jussieu, Paris Cedex 05, France.\\
e-mail: melleray@math.jussieu.fr
\end{document}